\documentclass[11pt]{article}
\usepackage{amsmath,amsfonts}
\usepackage{epsfig}

\def\ds{\displaystyle}
\makeatletter      
\def\@cite#1#2{[{%
  \m@th\upshape\mdseries{{\bfseries#1}\if@tempswa, #2\fi}}]}
\makeatother

\newcommand{\medio}{{\textstyle\frac{1}{2}}}
\newcommand{\qed}{\ifmmode
    \else\leavevmode\unskip\penalty9999 \hbox{}\nobreak\hfill\fi
	\quad\hbox{\qedsymbol}}
\newcommand{\qedsymbol}{\vbox{\hrule\hbox{\vrule height5.2pt
                              \hskip 5.2pt \vrule}\hrule}}
\newcommand{\regd}{$^{\mbox{\scriptsize\textregistered}}$}
\renewcommand{\tfrac}[2]{{\textstyle\frac{#1}{#2}}}

\begin{document}
\thispagestyle{empty}

\begin{center}
\Large
Location of Incenters and Fermat Points in Variable Triangles
\end{center}

\vspace*{0.2in}

\begin{flushright}
ANTHONY V\'ARILLY\\
Harvard University\\
Cambridge, MA 02138
\end{flushright}

\vskip .25in

\paragraph*{Introduction}
In 1765, Euler proved that several important centers of a triangle are
collinear; the line containing these points is named after him. The
incenter $I$ of a triangle, however, is generally not on this line.
Less than twenty years ago, Andrew P. Guinand discovered that although
$I$ is not necessarily on the Euler line, it is always fairly close to
it. Guinand's theorem~\cite{Guinand} states that for any
non-equilateral triangle, the incenter lies inside, and the excenters
lie outside, the circle whose diameter joins the centroid to the
orthocenter; henceforth the orthocentroidal circle. Furthermore,
Guinand constructed a family of locus curves for~$I$ which cover the
interior of this circle twice, showing that there are no other
restrictions on the positions of the incenter with respect to the
Euler line.

Here we show that the Fermat point also lies inside the
orthocentroidal circle; this suggests that the neighborhood of the
Euler line may harbor more secrets than was previously known. We also
construct a simpler family of curves for~$I$, covering the interior
once only except for the nine-point center~$N$, which corresponds to
the limiting case of an equilateral triangle.

Triangle geometry is often believed to be exhausted, although both
Davis~\cite{Davis} and Oldknow~\cite{Oldknow} have expressed the hope
that the use of computers may revive it. New results do appear
occasionally, such as Eppstein's recent construction of two new
triangle centers~\cite{Eppstein}. This article establishes some
relations among special points of the triangle, which were indeed
found by using computer software.

\paragraph*{The locus of the incenter}
To any plane triangle $ABC$ there are associated several special
points, called \emph{centers}. A few of these, in the standard
notation~\cite{Court}, are: the centroid~$G$, 
where the medians intersect;
the orthocenter~$H$, where the altitudes meet; 
the centers of the inscribed and circumscribed circles, called the
incenter~$I$ and the circumcenter~$O$; and
the nine-point center~$N$, half-way between~$O$ and~$H$. 
The
radii of the circumcircle and the incircle are called $R$ and~$r$,
respectively. 

If equilateral triangles $BPC$, $AQC$ and $ARB$ are
constructed externally on the sides of the triangle $ABC$, then the
lines $AP$, $BQ$ and $CR$ are concurrent and meet at the Fermat
point~$T$. This point minimizes the distance $TA + TB + TC$ for
triangles whose largest angle is $\leq 120^\circ$~\cite{Coxeter}.

The points $O$, $G$, $N$, and $H$ lie (in that order) on a line called the Euler
line, and $OG : GN : NH = 2 : 1 : 3$. They are distinct unless $ABC$
is equilateral. The circle whose diameter is $GH$ is called the
\emph{orthocentroidal circle}.

Guinand noticed that Euler's relation $OI^2 = R(R - 2r)$
\cite[p.~85]{Court} and Feuerbach's theorem $IN = \medio R - r$
\cite[p.~105]{Court} together imply that
$$
 OI^2 - 4\,IN^2 = R(R - 2r) - (R - 2r)^2 = 2r(R - 2r)
 = \frac{2r}{R} \, OI^2 > 0.
$$
Therefore, $OI > 2\,IN$. The locus of points $P$ for which
$OP = 2\,PN$ is a circle of Apollonius; since $OG = 2\,GN$ and
$OH = 2\,HN$, this is the orthocentroidal circle. The inequality
$OI > 2\,IN$ shows that $I$ lies in the interior of the
circle~\cite{Guinand}.

Guinand also showed that the angle cosines $\cos A$, $\cos B$,
$\cos C$ of the triangle satisfy the following cubic equation:
\begin{equation}
 \rho^4 (1 - 2x)^3 + 8\rho^2 \sigma^2 x(3 - 2x)
  - 16\sigma^4 x - 4\sigma^2 \kappa^2 (1 - x) = 0,
 \label{eq:Guinand}
\end{equation}
where $OI = \rho$, $IN = \sigma$ and $OH = \kappa$.  We exploit this
relationship below.

The relation $OI > 2\,IN$ can be observed on a computer with the
software \emph{The Geometer's Sketchpad}\regd, that allows tracking of
relative positions of objects as \emph{one} of them is moved around
the screen. Let us fix the Euler line by using a Cartesian coordinate
system with $O$ at the origin and $H$ at $(3,0)$. Consequently,
$G = (1,0)$ and $N = (1.5,0)$. To construct a triangle with this Euler
line, we first describe the circumcircle $\odot(O,R)$, centered at~$O$
with radius $R > 1$ ---in order that $G$ lie in the interior--- and
choose a point $A$ on this circle. If $AA'$ is the median passing
through $A$, we can determine $A'$ from the relation
$AG : GA' = 2 : 1$; then $BC$ is the chord of the circumcircle that is
bisected perpendicularly by the ray $OA'$.

\begin{figure}[ht]
  \begin{center}
    \epsfig{file=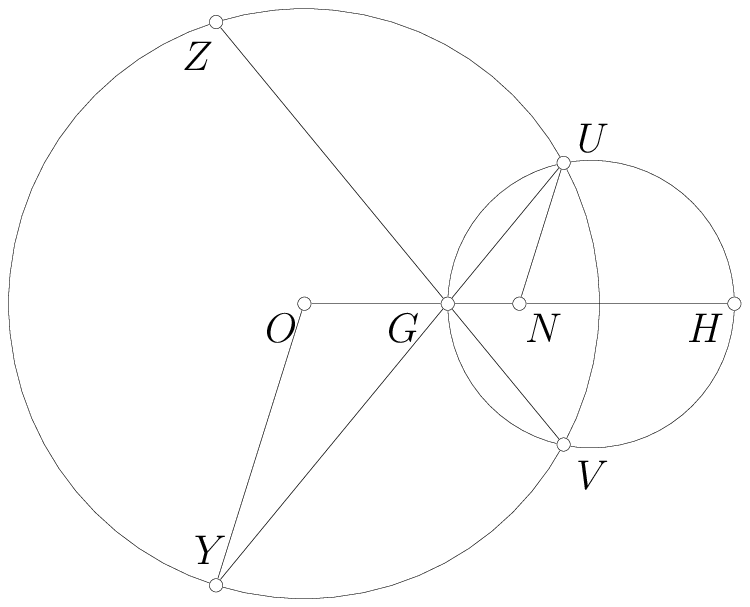}
    \caption{The forbidden arc $ZY$}
    \label{fg:forbidarc}
  \end{center}
\end{figure}


It is not always possible to construct $ABC$ given a fixed Euler line
and a circumradius $R$. If $1 < R < 3$ then there is an arc on
$\odot(O,R)$ on which an arbitrary point $A$ yields an $A'$ outside
the circumcircle, which is absurd. If $U$ and $V$ are the
intersections of $\odot(O,R)$ with the orthocentroidal circle, and
if $UY$ and $VZ$ are the chords of $\odot(O,R)$ passing through $G$,
then $A$ cannot lie on the arc $ZY$ of $\odot(O,R)$, opposite the
orthocentroidal circle (Figure~\ref{fg:forbidarc}). $OYG$ and $NUG$
are similar triangles. Indeed, $YO = UO = 2\,UN$ and $OG = 2\,GN$, so
$YG : GU = 2 : 1$; in the same way, $ZG : GV = 2 : 1$. Hence, if
$A = Y$ or $A = Z$, then $A' = U$ or $A' = V$ respectively, and $ABC$
degenerates into the chord $UY$ or $VZ$. If $A$ lies on the arc $ZY$,
opposite the orthocentroidal circle, $A'$ will be outside $\odot(O,R)$.

{}Another formula of Guinand~\cite{Guinand},
\begin{equation}
 OH^2 = R^2(1 - 8 \cos A \cos B \cos C),
 \label{eq:OHsq}
\end{equation}
shows that the forbidden arc appears if and only if $\cos A \cos B \cos C < 0$,
that is, if and only if $1 < R < 3 = OH$. The triangle $ABC$ is
obtuse-angled in this case.

Once $A$ is chosen and the triangle constructed, we can find $I$ by
drawing the angle bisectors of $ABC$ and marking their intersection.
\emph{The Geometer's Sketchpad} will do this automatically, and we can
also ask it to draw the locus which $I$ traces as we move the point
$A$ around the circumcircle. The idea of parameterizing the loci with
$R$, instead of an angle of the triangle (as Guinand did) was inspired
by the drawing tools available in \emph{The Geometer's Sketchpad}.
This locus turns out to be a quartic curve, as follows.

{\sc Proposition 1.}
{\it In the coordinate system described above, the incenter of $ABC$ is on the curve}
\begin{equation}
 (x^2 + y^2)^2 = R^2\,[(2x - 3)^2 + 4y^2].
 \label{eq:Tony}
\end{equation}

\textit{Proof}. From equation (\ref{eq:OHsq}) we see that, once we fix
$R$, the product $\cos A \cos B \cos C$ is fixed. Using Vi\`ete's
formulas, this product is obtained from the constant term
of Guinand's equation~(\ref{eq:Guinand}):
$$
 \cos A \cos B \cos C
  = \frac{1}{8} \left( 1 - \frac{4\sigma^2\kappa^2}{\rho^4} \right),
$$
so that
\begin{equation}
 \rho^4 (1 - 8 \cos A \cos B \cos C) = 4\sigma^2 \, OH^2.
 \label{eq:Viete}
\end{equation}
Now, $I$ is a point of intersection of the two circles
$$
 \rho^2 = x^2 + y^2   \quad\mbox{and}\quad
 \sigma^2 = (x - \tfrac{3}{2})^2 + y^2.
$$
We get~(\ref{eq:Tony}) by substituting these and (\ref{eq:OHsq}) in
(\ref{eq:Viete}), and dividing through by the common factor
$1 - 8 \cos A \cos B \cos C$, which is positive since $O \neq H$
(the equilateral case is excluded).   \qed

\vspace{6pt}

\begin{figure}[ht]
  \begin{center}
    \epsfig{file=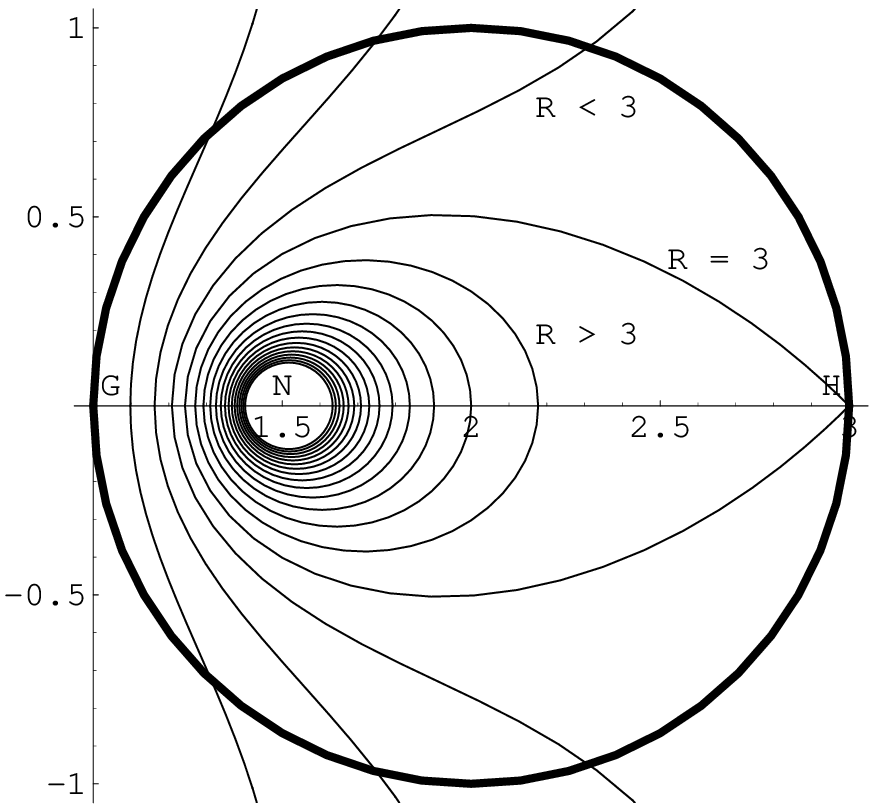}
    \caption{Locus curves for $I$ with fixed $R$}
    \label{fg:curves}
  \end{center}
\end{figure}


Figure~\ref{fg:curves} is a \emph{Mathematica} plot of the curves
(\ref{eq:Tony}) with the orthocentroidal circle.

Every point on a locus inside the orthocentroidal circle is the
incenter of a triangle. When $R > 3$, the locus is a lobe entirely
inside the circle; if the point $A$ travels around the circumcircle
once, then $I$ travels around the lobe three times, since $A$ will
pass through the three vertices of each triangle with
circumradius~$R$. When $1 < R < 3$, the interior portion
of~(\ref{eq:Tony}) isshaped like a bell (Figure~\ref{fg:curves}). Let $A$
travel along the allowable arc from $Y$ to $Z$, passing through $V$
and $U$; then $I$ travels along the bell from $U$ to $V$, back from
$V$ to $U$, and then from $U$ to $V$ again. While $A$ moves from $V$
to~$U$, the orientation of the triangle $ABC$ is reversed. When
$R = 3$, the locus closes at $H$ and one vertex $B$ or $C$ also
coincides with $H$; the triangle is right-angled. If $A$ moves once
around the circumcircle, starting and ending at $H$, $I$ travels twice
around the lobe.

{\sc Proposition 2.}
{\it The curves (\ref{eq:Tony}), for different values of $R$, do not
cut each other inside the orthocentroidal circle, and they fill the
interior except for the point~$N$.}

\textit{Proof}. Let $(a,b)$ be inside the orthocentroidal circle, that
is,
\begin{equation}
 (a - 2)^2 + b^2 < 1.
 \label{eq:GuinandDisc}
\end{equation}
If $(a,b)$ also lies on one of the curves (\ref{eq:Tony}), then
$$
R = \sqrt{\frac{(a^2 + b^2)^2}{(2a - 3)^2 + 4b^2}}.
$$
There is only one positive value of $R$ and thus \emph{at most} one
curve of the type~(\ref{eq:Tony}) on which $(a,b)$ can lie. Now we
show $(a,b)$ lies on \emph{at least} one curve of
type~(\ref{eq:Tony}); to do that, we need to show that given
(\ref{eq:GuinandDisc}), $R > 1$. We need only prove
\begin{equation}
 (a^2 + b^2)^2 > (2a - 3)^2 + 4b^2.
 \label{eq:Rgeone}
\end{equation}
Indeed, $(2a - 3)^2 + 4b^2 = 0$ only if $(a,b) = (\frac{3}{2},0)$;
this point is $N$. It cannot lie on a locus (\ref{eq:Tony}), in fact
it corresponds to the limiting case of an equilateral triangle as
$R \to \infty$.

The inequality (\ref{eq:GuinandDisc}) can be restated as
\begin{equation}
 a^2 + b^2 < 4a - 3,
 \label{eq:reGuinandDisc}
\end{equation}
and (\ref{eq:Rgeone}) as
\begin{equation}
 {(a^2 + b^2)}^2 - 4a^2 + 3(4a - 3) - 4b^2 > 0.
 \label{eq:reRgeone}
\end{equation}
{}From (\ref{eq:reGuinandDisc}) it follows that
\begin{eqnarray*}
 (a^2 + b^2)^2 - 4a^2 + 3(4a - 3) - 4b^2
 & > & (a^2 + b^2)^2 - 4a^2 + 3(a^2 + b^2) - 4b^2  \\
 & = & (a^2 + b^2)(a^2 + b^2 - 1).
\end{eqnarray*}
But $a^2 + b^2 > 1$ since $(a,b)$ is inside the orthocentroidal
circle; therefore, (\ref{eq:Rgeone}) is true. \qed

\paragraph*{The whereabouts of the Fermat point}
The same set-up, a variable triangle with fixed Euler line and
circumcircle, allows us to examine the loci of other triangle centers.
Further experimentation with \emph{The Geometer's Sketchpad} suggests
that the Fermat point $T$ also lies inside the orthocentroidal circle
in all cases.

\begin{figure}[ht]
  \begin{center}
    \epsfig{file=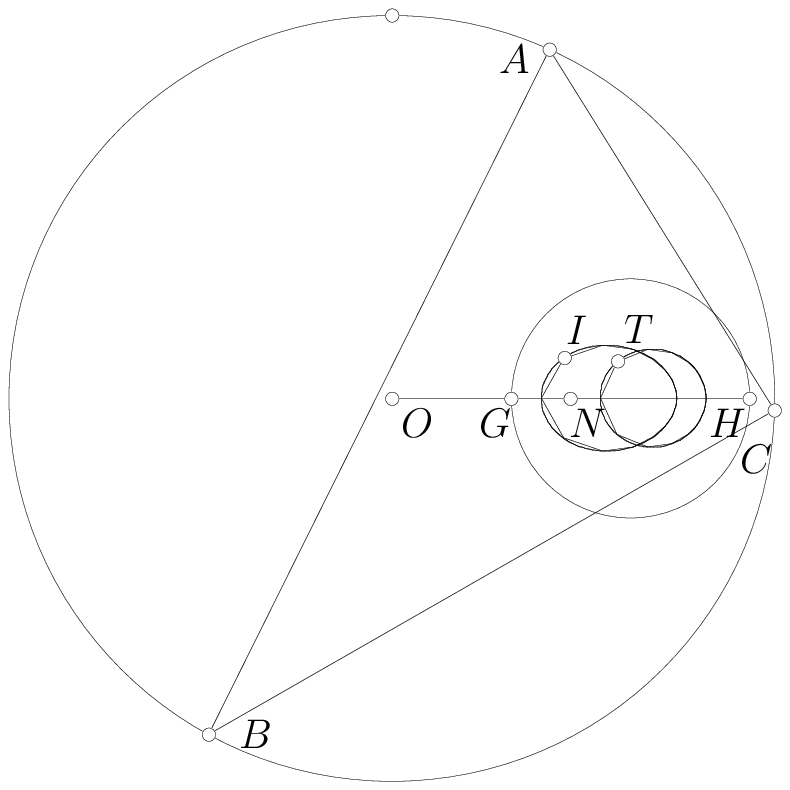}
    \caption{Locus of the Fermat point $T$}
    \label{fg:Ferlocus}
  \end{center}
\end{figure}


{\sc Theorem 1.}
{\it The Fermat point of any non-equilateral triangle lies inside the
orthocentroidal circle.}

\begin{figure}[ht]
  \begin{center}
    \epsfig{file=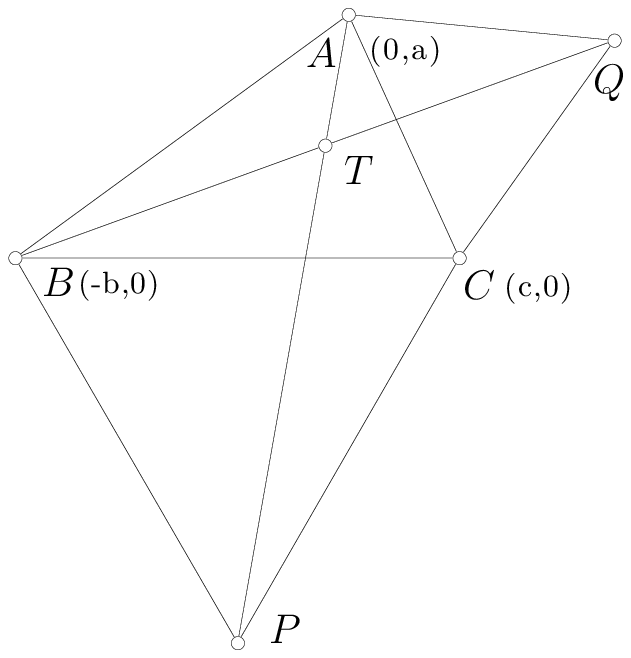}
    \caption{Finding the Fermat point}
    \label{fg:Fermatpt}
  \end{center}
\end{figure}

\textit{Proof}. Given a triangle $ABC$ with the largest angle at $A$,
we can take a coordinate system with $BC$ as the $x$-axis and $A$ on
the $y$-axis. Then $A = (0,a)$, $B = (-b,0)$, $C = (c,0)$ where $a$,
$b$ and $c$ are all positive. Let $BPC$ and $AQC$ be the equilateral
triangles constructed externally over $BC$ and $AC$ respectively
(Figure~\ref{fg:Fermatpt}). Then
$P = (\medio(c-b), -\medio\sqrt{3}(b+c))$ and
$Q = (\medio(\sqrt{3}\,a + c), \medio(a + \sqrt{3}\,c))$.

The coordinates of $T$ can be found by writing down the equations for
the lines $AP$ and $BQ$ and solving them simultaneously. After a
little work, we get
$\ds T = \left(\frac{u}{d}, \frac{v}{d}\right)$, where
\begin{eqnarray}
 u & = & (\sqrt{3}\,bc - \sqrt{3}\,a^2 - ac - ab)(b - c),
 \nonumber \\
 v & = & (a^2 + \sqrt{3}\,ab + \sqrt{3}\,ac + 3bc)(b + c),
 \nonumber \\
 d & = & 2\sqrt{3}(a^2 + b^2 + c^2) + 6ac + 6ab + 2\sqrt{3}\,bc.
 \label{eq:Fermatcoords}
\end{eqnarray}

The perpendicular bisectors of $BC$ and $AC$ intersect at the
circumcenter $O = (\medio(c-b), (a^2 - bc)/2a)$. The nine-point
center $N$ is the  circumcenter of the triangle whose vertices are the
midpoints of the sides of $ABC$; we can deduce that
$N = (\frac{1}{4}(c-b), (a^2 + bc)/4a)$.

To show that $T$ lies inside the orthocentroidal circle, we must show
that $OT > 2\,NT$, or
\begin{equation}
 OT^2 > 4\,NT^2.
 \label{eq:Fermat}
\end{equation}
In coordinates, this inequality takes the equivalent form:
$$
 \left[\frac{u}{d} - \left( \frac{c-b}{2}\right) \right]^2
   + \left[\frac{v}{d} - \left( \frac{a^2 - bc}{2a} \right) \right]^2
 > 4\left[\frac{u}{d} - \left( \frac{c-b}{4} \right) \right]^2
   + 4\left[\frac{v}{d} - \left( \frac{a^2 + bc}{4a} \right) \right]^2,
$$
or, multiplying by $(2ad)^2$,
$$
 [2au - ad(c - b)]^2 + [2av - d(a^2 - bc)]^2
 > [4au - ad(c - b)]^2 + [4av - d(a^2 + bc)]^2.
$$
After expanding and canceling terms, this simplifies to
$$
 4a^2du(c - b) + 4adv[(2a^2 + 2bc) - (a^2 - bc)] - 4a^2d^2bc
 > 12a^2u^2 + 12a^2v^2,
$$
or better,
\begin{equation}
 adu(c - b) + dv(a^2 + 3bc) - abcd^2  - 3au^2 - 3av^2 > 0.
 \label{eq:FermatLHS}
\end{equation}

One way to verify this inequality is to feed the equations
(\ref{eq:Fermatcoords}) into \emph{Mathematica}, which expands and
factors the left hand side of (\ref{eq:FermatLHS}) as
$$
 2(b + c)(\sqrt{3}\,a^2 + \sqrt{3}\,b^2 + \sqrt{3}\,c^2 + \sqrt{3}\,bc
           + 3ab + 3ac)(a^4 + a^2b^2 - 8a^2bc + a^2c^2 + 9b^2c^2).
$$
The first three factors are positive. The fourth factor can be
expressed as the sum of two squares,
$$
 (a^2 - 3bc)^2 + a^2(b - c)^2,
$$
and could be zero only if $a^2 = 3bc$ and $b = c$, so that
$a = \sqrt{3}\,b$. This gives an equilateral triangle with side $2b$.
Since the equilateral case is excluded, all the factors are positive,
which shows that (\ref{eq:FermatLHS}) is true, and therefore
(\ref{eq:Fermat}) holds.   \qed

\vspace{6pt}

Varying the circumradius $R$ and the position of the vertex $A$ with
\emph{The Geometer's Sketchpad} reveals a striking parallel between
the behavior of the loci of $T$ and those of $I$. It appears that the
loci of $T$ also foliate the orthocentroidal disc, never cutting each
other, in a similar manner to the loci of $I$
(Figure~\ref{fg:curves}). The locus of $T$ becomes a lobe when
$R = 3$, as is the case with the locus of $I$. Furthermore, the loci
of $T$ close in on the center of the orthocentroidal circle as
$R \to \infty$, just as the loci of $I$ close in on $N$
(Figure~\ref{fg:curves}).

It is difficult to prove these assertions with the same tools used to
characterize the loci of $I$, because we lack an equation analogous
to~(\ref{eq:Guinand}) involving $T$ instead of~$I$. A quick
calculation for non-equilateral isosceles triangles, however, shows
that $T$ can be anywhere on the segment $GH$ except for its midpoint.
This is consistent with the observation that the loci of $T$ close in
on the center of the orthocentroidal circle.

Consider a system of coordinates like those of
Figure~\ref{fg:Fermatpt}. Let $b = c$ so that $ABC$ is an isosceles
triangle. In this case, by virtue of~(\ref{eq:Fermatcoords}),
$T = (0, b/\sqrt{3})$. For this choice of coordinates, $G = (0, a/3)$
and $H = (0, b^2/a)$. $T$ lies on the Euler line, which can be
parametrized by $(1 - t)G + tH$ for real~$t$. This requires that
$$
 (1 - t)\frac{a}{3} + t\,\frac{b^{2}}{a} = \frac{b}{\sqrt{3}},
$$
for some real $t$. Solving for $t$, this becomes
$$
 t = \frac{a^{2} - \sqrt{3}\,ab}{a^{2} - 3b^{2}}
   = \frac{a}{a + \sqrt{3}\,b},
$$
unless $a = \sqrt{3}\,b$. This case is excluded since $ABC$ is not
equilateral. Note that $t \to \frac{1}{2}$ as $a \to \sqrt{3}\,b$.
Thus, $t$ takes real values between $0$ and $1$ except for
$\frac{1}{2}$, so $T$ can be anywhere on the segment $GH$ except for
its midpoint.

\vspace{6pt}

{\small
\paragraph*{Acknowledgments}
I am grateful to Les Bryant and Mark Villarino for helpful
discussions, to Joseph C. V\'arilly for advice and \TeX{}nical help,
and to Julio Gonz\'alez Cabill\'on for good suggestions on the use of
computer software. I also thank the referees for helpful comments.
}

\end{document}